\documentclass[a4paper,10pt]{article}
\usepackage[utf8]{inputenc}
\usepackage[margin=1.3in]{geometry}

\usepackage{enumerate}

\usepackage[english]{babel}
\usepackage{fancyhdr}
\usepackage{amsmath}
\usepackage{comment}
\let\oldmarginpar\marginpar
\renewcommand\marginpar[1]{\-\oldmarginpar[\raggedleft\small\sf
#1] {\raggedright\small\sf #1}}
\usepackage{color}
\usepackage{amssymb}
\usepackage{url}
\usepackage{alltt}
\usepackage[T1]{fontenc}
\usepackage{latexsym}
\usepackage{t1enc}
\usepackage{amsfonts}
\usepackage{amssymb}
\usepackage{enumerate}
\usepackage{stackrel}
\usepackage{bm}
\usepackage{upref}
\usepackage{comment,ulem}
\usepackage{graphics,psfrag,rotating}
\usepackage{cite}
\usepackage{ntheorem}

\newcommand{\beq}{\begin{equation}}
\newcommand{\eeq}{\end{equation}}

\newcommand{\p}{p(\cdot)}
\newcommand{\q}{q(\cdot)}
\newcommand{\supp}{\mathop{\mathrm{supp}\,}\nolimits}

\newcommand{\ds}{\displaystyle}
\newcommand{\real}{\mathbb{R}}
\newcommand{\rn}{\mathbb{R}^n}
\newcommand{\nat}{\mathbb{N}}

\newcommand{\zn}{\mathbb{Z}^n}

\newcommand{\bit}{\begin{itemize}}
\newcommand{\eit}{\end{itemize}}

\newtheorem{definition}{Definition }[section]

\newtheorem{theorem}[definition]{Theorem}

\newtheorem{lemma}[definition]{Lemma}
\theoremstyle{definition}\newtheorem{remark}[definition]{Remark}

\def\ds{\displaystyle}

\def\supp{\mbox{\rm supp }}

\newcommand{\bpr}{\begin{proof}}
\newcommand{\epr}{\nopagebreak \hfill $\square$\\ \end{proof}}

\title{Intrinsic atomic characterization of 2-microlocal spaces with variable exponents on domains}
\author{Helena Gon\c{c}alves \and Henning Kempka\thanks{The authors were supported by the German science foundation (DFG) within the project KE 1847/1-2.}}

\begin{document}
\maketitle

\begin{abstract}
We provide an intrinsic atomic characterization for 2-microlocal Besov and Triebel-Lizorkin spaces with variable integrability on domains, $B_{p(\cdot),\q}^{\bm{w}}(\Omega)$ and $F_{p(\cdot),\q}^{\bm{w}}(\Omega)$, where $\Omega$ is a regular domain. We make use of the non-smooth atomic decomposition result obtained in \cite{GK15} for these spaces to get the main result.
\end{abstract}


\section{Introduction}

\indent

The 2-microlocal function spaces initially appeared in the book of Peetre \cite{Pe76} and have also been studied by Bony \cite{Bo86} in connection with pseudodifferential operators. Later on, they were investigated by Jaffard \cite{Jaf91} as well as Jaffard and Meyer \cite{JM96}. In \cite{LevySeu03} and \cite{LevySeu04}, Levy V\'ehel and Seuret showed that they are a useful tool to measure local regularity and to describe the oscillatory behavior of functions near singularities.

Spaces of variable integrability, also known as variable exponent function spaces $L_{\p}(\rn)$, can be traced back to Orlicz \cite{Or31} 1931, but the modern development started with the papers \cite{KR91} of Kov\'a\v{c}ik and R\'akosn\'ik as well as \cite{Die04} of Diening. The spaces $L_{\p}(\rn)$ have interesting applications in fluid dynamics, image processing, PDE and variational calculus, see the introduction of \cite{DHR09}. For an overview we refer to \cite{DHHR}.

The concept of function spaces with variable smoothness and the concept of variable integrability were firstly mixed up by Diening, H{\"a}st{\"o} and Roudenko in \cite{DHR09}. They defined Triebel-Lizorkin spaces $F^{s(\cdot)}_{\p, \q}(\rn)$ and from the trace theorem on $\mathbb{R}^{n-1}$ it became clear why it is natural to have all parameters variable. Due to
\[
\mathrm{Tr}\  F^{s(\cdot)}_{\p,\q}(\rn)={F}^{s(\cdot)-\frac{1}{\p}}_{\p,\p}(\mathbb{R}^{n-1}), \mbox{ with } s(\cdot) - \frac{1}{\p} > (n-1)\max\left(\frac{1}{\p}-1, 0\right),
\]
(\!\!\cite[Theorem~3.13]{DHR09}) we see the necessity of taking $s$ and $q$ variable if $p$ is not constant. A similar interplay between smoothness and integrability is known to be inherited by the Sobolev embedding. The Sobolev embedding could be transferred to the case of variable integrability by Leopold and Edmunds \& R\'akosn\'\i k in \cite{Leopold}, \cite{Ed1} and \cite{Ed2} and the corresponding result for Triebel-Lizorkin spaces with variable smoothness and integrability has been given by Vyb\'iral in \cite{Vyb09}. Moreover, Almeida and H{\"a}st{\"o} also introduced in \cite{AH10} Besov spaces $B_{\p,\q}^{s(\cdot)}(\rn)$ with all three indices variable and showed a Sobolev embedding for these spaces.

The scale we consider here - mixing admissible weight sequences with variable integrability - was introduced in \cite{Kem09, Kem10} and provides a unified approach that covers many spaces related with variable smoothness and generalized smoothness. Many results have been studied regarding these spaces, in particular the possibility of decomposing functions $f \in B_{\p,\q}^{\bm{w}}(\rn)$ or $F_{\p,\q}^{\bm{w}}(\rn)$ as linear combinations of smooth atoms, which are the building blocks for atomic decompositions. More recently, a more general decomposition for these spaces was obtained in \cite{GK15}, where the authors show that one can replace the usual atoms used in smooth atomic decompositions by more general ones, called non-smooth atoms. Those atoms are characterized by a relaxation on the smoothness assumptions and, nevertheless, one keeps all the crucial information compared to smooth atomic decompositions. We devote Section 3 to this topic.

Regarding intrinsic characterizations of function spaces on domains, in \cite{TW96} Triebel and Winkelvoß suggested the use of these non-smooth atoms as a tool to define classical Besov and Triebel-Lizorkin spaces $B^s_{p,q}(\Omega)$ and $F^s_{p,q}(\Omega)$ on a class of (non-smooth) domains. Also Rychkov in \cite{Ry98} gave an intrinsic characterization for the same scale of spaces, but considering smooth domains. More recently, Tyulenev in \cite{Tyu16} studied Besov-type spaces of variable smoothness on rough domains, namely bounded Lipschitz domains in $\rn$, epigraph of Lipschitz functions or $(\epsilon, \delta)$-domains. Concerning 2-microlocal Besov and Triebel-Lizorkin spaces with variable exponents, Kempka presented recently in \cite{Kem16} two different intrinsic characterizations of these spaces using local means and the Peetre maximal operator, on special Lipschitz domains.

Since a non-smooth atomic characterization for the scale of 2-microlocal Besov and Triebel-Lizorkin spaces $B_{\p,\q}^{\bm{w}}(\rn)$ and $F_{\p,\q}^{\bm{w}}(\rn)$ was already obtained, our aim is to get an intrinsic characterization of these spaces for more general domains, as considered in \cite{TW96}. We deal with this problem in Section 4, where we study spaces on the scale of \textit{regular domains}. We wish to emphasize that this class of domains includes not only bounded connected Lipschitz domains but also special Lipschitz domains and $(\epsilon, \delta)$-domains.

\section{Notation and definitions}

\indent

We shall adopt the following general notation: $\nat$ denotes the set of all natural numbers, $\nat_0=\mathbb N\cup\{0\}$, $\mathbb{Z}$ denotes the set of integers, $\rn$ for $n\in\nat$ denotes the $n$-dimensional real Euclidean space with $|x|$, for $x\in\rn$, denoting the Euclidean norm of $x$. For a real number $a$, let $a_+:=\max(a,0)$.

If $s \in \mathbb{R}$, then there are uniquely determined $\lfloor s\rfloor^- \in \mathbb{Z}$ and $\{s\}^+ \in (0,1]$ with $s = \lfloor s\rfloor^- + \{s\}^+$.

\begin{definition}\label{holder-space}
  Let $s >0$. Then the H{\"o}lder space with index $s$ is defined as
  \beq
    \mathcal{C}^s(\rn) = \Big\{ f \in C^{\lfloor s\rfloor^-} (\rn): \| f \mid \mathcal{C}^s (\rn)\| < \infty \Big\}, \nonumber
  \eeq
  with
  \beq
    \| f \mid \mathcal{C}^s (\rn)\|:= \sum_{|\alpha|\leq \lfloor s\rfloor^-} \sup_{x \in \rn} |D^{\alpha}f(x)| + \sum_{|\alpha| = \lfloor s\rfloor^-} \sup_{x,y \in \rn, x\neq y} \frac{|f(x)-f(y)|}{|x-y|^{\{s\}^+}}. \nonumber
  \eeq
\end{definition}
If $s=0$, then we set $\mathcal{C}^0(\rn)=L_{\infty}(\rn)$.

\indent

For $q\in (0,\infty]$, $\ell_q$ stands for the linear space of all complex sequences $f=(f_j)_{j\in \nat_0}$ endowed with the quasi-norm
$$
\Vert f \mid \ell_q\Vert = \Big(\sum_{j=0}^{\infty} |f_j|^q \Big)^{1/q},
$$
with the usual modification if $q=\infty$. By $c$, $c_1$, $c_2$, etc. we denote positive constants independent of appropriate quantities. For two non-negative expressions ({\it i.e.}, functions or functionals) ${\mathcal  A}$, ${\mathcal  B}$, the symbol ${\mathcal A}\lesssim {\mathcal  B}$ (or ${\mathcal A}\gtrsim {\mathcal  B}$) means that $ {\mathcal A}\leq c\, {\mathcal  B}$ (or $c\,{\mathcal A}\geq {\mathcal B}$), for some $c>0$. If ${\mathcal  A}\lesssim {\mathcal  B}$ and ${\mathcal A}\gtrsim{\mathcal  B}$, we write ${\mathcal  A}\sim {\mathcal B}$ and say that ${\mathcal  A}$ and ${\mathcal  B}$ are equivalent.

In order to define 2-microlocal Besov and Triebel-Lizorkin spaces with variable integrability, we start by recalling the definition of admissible weight sequences. We follow \cite{Kem10}.

\begin{definition}\label{def-ad-weight} Let $\alpha\geq 0$ and $\alpha_1,\alpha_2\in\real$ with $\alpha_1\leq \alpha_2$. A
sequence of non-negative measurable functions in $\rn$ $\bm{w}=(w_j)_{j\in\nat_0}$ belongs to the class $\mathcal{W}^{\alpha}_{\alpha_1,\alpha_2}(\rn)$ if the following conditions are satisfied:
  \begin{list}{}{\labelwidth1.7em\leftmargin2.3em}
    \item[{\hfill (i)\hfill}] There exists a constant $c>0$ such that
      $$
	0<w_j(x)\leq c\,w_j(y)\,(1+2^j|x-y|)^{\alpha}\quad \mbox{for all} \;\, j\in\nat_0 \;\; \mbox{and all} \;\, x,y\in\rn.
      $$
    \item[{\hfill (ii)\hfill}] For all $j\in\nat_0$ it holds
      $$
	2^{\alpha_1}\,w_j(x)\leq w_{j+1}(x)\leq 2^{\alpha_2}\,w_j(x) \quad \mbox{for all}\;\, x \in \rn.
      $$
  \end{list}
Such a system $(w_j)_{j\in\nat_0}\in\mathcal{W}^{\alpha}_{\alpha_1,\alpha_2}(\rn)$ is called admissible weight sequence.\\
\end{definition}

Properties of admissible weights may be found in \cite[Remark~2.4]{Kem08}.

Before introducing the function spaces under consideration we still need to recall some notation. By $\mathcal{S}(\rn)$ we denote the Schwartz space of all complex-valued rapidly decreasing infinitely differentiable functions on $\rn$ and by $\mathcal{S}'(\rn)$ the dual space of all tempered distributions on $\rn$. For $f\in \mathcal{S}'(\rn)$ we denote by $\widehat{f}$ the Fourier transform of $f$ and by $f^{\vee}$ the inverse Fourier transform of $f$.

Let $\varphi_0\in\mathcal{S}(\rn)$ be such that
\begin{equation}  \label{phi}
  \varphi_0(x)=1 \quad \mbox{if}\quad |x|\leq 1 \quad \mbox{and} \quad \supp \varphi_0 \subset \{x\in\rn: |x|\leq 2\}.
\end{equation}
Now define $\varphi(x):=\varphi_0(x)-\varphi_0(2x)$ and set $\varphi_j(x):=\varphi(2^{-j}x)$ for all $j\in\nat$. Then the sequence $(\varphi_j)_{j\in\nat_0}$ forms a smooth dyadic partition of unity.

By $\mathcal{P}(\rn)$ we denote the class of exponents, which are measurable functions $p:\rn\rightarrow (c,\infty]$ for some $c>0$. Let $p \in \mathcal{P}(\rn)$. Then, $p^+:=\mbox{ess-sup}_{x\in \rn}p(x)$, $p^-:=\mbox{ess-inf}_{x\in \rn}p(x)$ and $L_{p(\cdot)}(\rn)$ is the variable exponent Lebesgue space, which consists of all measurable functions $f$ such that for some $\lambda>0$ the (quasi-)modular $\varrho_{L_{p(\cdot)}(\rn)}(f/\lambda)$ is finite, where
$$
\displaystyle \varrho_{L_{p(\cdot)}(\rn)}(f):=\int_{\rn_0 } |f(x)|^{p(x)}\, dx + \mbox{ess-sup}_{x\in \rn_{\infty}} |f(x)|.
$$
Here $\rn_{\infty}$ denotes the subset of $\rn$ where $p(x)=\infty$ and $\rn_0= \rn \setminus \rn_{\infty}$. The Luxemburg quasi-norm (norm if $p(x)\geq1$) of a function $f\in L_{p(\cdot)}(\rn)$ is given by
$$
\|f \mid L_{p(\cdot)}(\rn)\|:=\inf\left\{\lambda>0:\varrho_{L_{p(\cdot)}(\rn)}\left(\frac{f}{\lambda}\right)\leq 1 \right\}.
$$

In order to define the mixed spaces $\ell_{\q}(L_{\p})$, we need to define another (quasi-)modular. For $p,q \in \mathcal{P}(\rn)$ and a sequence $(f_{\nu})_{\nu \in \nat_0}$ of complex-valued Lebesgue measurable functions on $\rn$, we define
\beq \label{modular_mixed}
  \varrho_{\ell_{\q}(L_{\p})}(f_{\nu}) = \sum_{\nu=0}^{\infty} \inf\left\{\lambda_{\nu}>0 : \varrho_{\p}\left(\frac{f_{\nu}}{\lambda_{\nu}^{1/\q}} \right)\leq 1 \right\}.
\eeq
If $q^+ <\infty$, then we can replace \eqref{modular_mixed} by the simpler expression
\beq \label{modular_mixed_norm}
  \varrho_{\ell_{\q}(L_{\p})}(f_{\nu}) = \sum_{\nu=0}^{\infty} \Big\| |f_{\nu}|^{\q} \mid L_{\frac{\p}{\q}} \Big\|.
\eeq
The (quasi-)norm in the $\ell_{\q}(L_{\p})$ spaces is defined as usual by
\beq
  \| f_{\nu} \mid \ell_{\q}(L_{\p}(\rn)) \| = \inf \left\{ \mu>0 : \varrho_{\ell_{\q}(L_{\p})}\left(\frac{f_{\nu}}{\mu}\right) \leq 1\right\}.
\eeq

For the sake of completeness, we state also the definition of the space $L_{\p}(\ell_{\q})$. At first, one just takes the norm $\ell_{\q}$ of $(f_{\nu}(x))_{\nu \in \nat_0}$ for every $x \in \rn$ and then the $L_{\p}$-norm with respect to $x \in \rn$, i.e.
\beq
  \| f_{\nu} \mid L_{\p}(\ell_{\q}(\rn)) \| = \left\| \left(\sum_{\nu=0}^{\infty} |f_{\nu}(x)|^{q(x)} \right)^{1/q(x)} \mid L_{\p}(\rn)\right\|. \nonumber
\eeq \\


\begin{definition}\label{log-holder}
  Let $g\in C(\rn)$. We say that $g$ is locally log-H\"older continuous, abbreviated $g\in C_{loc}^{\log}(\rn)$, if there exists $c_{\log}(g)\geq0$ such that
  \beq \label{loc-log-holder}
      |g(x)-g(y)|\leq \frac{c_{\log}(g)}{\log(e+1/|x-y|)} \quad \text{for all}\;\; x,y\in\rn.
   \eeq
  We say that $g$ is globally log-H\"older continuous, abbreviated $g\in C^{\log}(\rn)$, if $g$ is locally log-H\"older continuous and there exists $g_{\infty}\in\real$ and $c_{\log}\geq0$ such that
    \beq
      |g(x)-g_{\infty}|\leq \frac{c_{\log}}{\log(e+|x|)} \quad \text{for all}\;\; x\in\rn.
    \eeq
\end{definition}

We use the notation $p\in \mathcal{P}^{\log}(\rn)$ if $p\in \mathcal{P}(\rn)$ and $1/p\in C^{\log}(\rn)$. 

The definitions of the spaces below were given in \cite{KemV12}.
\begin{definition}  \label{defspaces}
Let $(\varphi_j)_{j\in\nat_0}$ be a partition of unity as above, $\bm{w}=(w_j)_{j\in\nat_0}\in\mathcal{W}^{\alpha}_{\alpha_1,\alpha_2}(\rn)$ and $p,q \in \mathcal{P}^{\log}(\rn)$.
\begin{list}{}{\labelwidth1.3em\leftmargin2em}
  \item[{\upshape (i)\hfill}] The space $B_{p(\cdot),\q}^{\bm{w}}(\rn)$ is defined as the collection of all $f\in \mathcal{S}'(\rn)$ such that
    \begin{align}
      \|f\mid B_{p(\cdot),\q}^{\bm{w}}(\rn)\|_{\varphi} &:= \Vert (w_j\,(\varphi_j \widehat{f})^{\vee} )_{j\in\nat_0}\mid \ell_{\q}(L_{p(\cdot)}(\rn))\Vert\notag
    \end{align}
    is finite.
  \item[{\upshape (ii)\hfill}] If $p^{+}, q^+< \infty$, then the space $F_{p(\cdot),q(\cdot)}^{\bm{w}}(\rn)$ is defined as the collection of all $f\in \mathcal{S}'(\rn)$ such that
    \begin{align}
      \|f \mid F^{\bm{w}}_{\p,\q}(\rn) \|_{\varphi} &:= \Vert (w_j\,(\varphi_j \widehat{f})^{\vee} )_{j\in\nat_0}\mid L_{p(\cdot)}(\ell_{q(\cdot)} (\rn))\Vert\notag
    \end{align}
    is finite.
\end{list}
\end{definition}

\begin{remark}
    These spaces include very well-known spaces. For $p=$const and $w_j(x)=2^{js}$ we get back to the classical Besov and Triebel-Lizorkin spaces $B^{s}_{p,q}(\rn)$ and $F^{s}_{p,q}(\rn)$.

    Also the spaces of generalized smoothness are contained in this approach (see \cite{FL06}, \cite{Mou01}) by taking
    $$w_j(x) = 2^{js} \Psi(2^{-j}), \qquad \mbox{ or more general } w_j(x)=\sigma_j.$$
    Here, $\{\sigma_j\}_{j \in \nat_0}$ is an admissible sequence, which means that there exist $d_0, d_1>0$ with $d_0 \sigma_j \leq \sigma_{j+1} \leq d_1 \sigma_j$ and $\Psi$ is a slowly varying function.

    Moreover, these 2-microlocal spaces also cover the spaces of variable smoothness and integrability $B^{s(\cdot)}_{\p, \q}$ and $F^{s(\cdot)}_{\p, \q}$, introduced in \cite{DHR09} and \cite{AH10}. If $s \in C_{loc}^{\log}(\rn)$ (which is the standard condition on $s(\cdot)$), then $\bm{w}= (w_j(x))_{j \in \nat_0}= (2^{js(x)})_{j \in \nat_0}$ belongs to $\mathcal{W}^{\alpha}_{\alpha_1,\alpha_2}(\rn)$ with $\alpha_1=s^-$, $\alpha_2=s^+$ and $\alpha=c_{log}(s)$, where $c_{log}(s)$ is the constant for $s(\cdot)$ from \eqref{loc-log-holder}.
\end{remark}

For $p, q\in\mathcal{P}(\rn)$, we put
$$
\sigma_{p}:=n\left(\frac{1}{p^-}-1\right)_+ \quad \mbox{and}\quad \sigma_{p,q}:=n\left(\frac{1}{\min(1,p^-, q^-)}-1\right).
$$

\section{Non-smooth atomic characterization}

\indent

In this section we present a non-smooth atomic decomposition result for 2-microlocal Besov and Triebel-Lizorkin spaces $B_{\p,\q}^{\bm{w}}(\rn)$ and $F_{\p,\q}^{\bm{w}}(\rn)$, proved in \cite{GK15}.

At first, we shall introduce some notation. Let $\zn$ stand for the lattice of all points in $\rn$ with integer-valued components. Let $b>0$ be given, $\nu \in \nat_0$ and $m=(m_1,\dots,m_n) \in \zn$. Then $Q_{\nu,m}$ denotes a cube in $\rn$ with sides parallel to the axes of coordinates, centered at $x^{\nu,m} \in \rn$ with
\beq \label{defcenter}
    | x^{\nu,m} - 2^{- \nu}m| \leq b \, 2^{-\nu}
\eeq
and with side length $2^{-\nu}$. If $Q$ is a cube in $\rn$ and $r>0$, then $r\,Q$ is the cube in $\rn$ concentric with $Q$ and with side length $r$ times the side length of $Q$. By $\chi_{\nu,m}$ we denote the characteristic function of the cube  $Q_{\nu,m}$. In the sequel, we always implicitly assume that $d>0$ is chosen in dependence on $b$ such that for all choices of $\nu \in \nat_0$ and all choices of $x^{\nu,m}$ in \eqref{defcenter}
\beq \label{d-param}
    \bigcup_{m \in \zn} d\, Q_{\nu,m} = \rn.
\eeq

\begin{definition}\label{defatoms}
  Let $K, L\geq 0$, $d>1$ and $c>0$. A function $a: \rn \rightarrow \mathbb{R}$ is called a non-smooth $[K, L]$-atom centered at $Q_{\nu,m}$, for all $\nu \in \nat_0$ and $m \in \zn$, if
  \begin{list}{}{\labelwidth1.3em\leftmargin2em}
    \item[{\upshape (i)\hfill}] $\supp a \subset d\,Q_{\nu,m}$,
    \item[{\upshape (ii)\hfill}] $ \|a(2^{-\nu} \cdot) \mid \mathcal{C}^K (\rn) \| \leq c\,$,
    \item[{\upshape (iii)\hfill}] and for every $\psi \in \mathcal{C}^L(\rn)$ it holds $$\Big| \ds \int_{d\,Q_{\nu,m}} \psi(x) a(x) dx\Big| \leq c \, 2^{-\nu(L + n)}\| \psi \mid \mathcal{C}^L(\rn)\|.$$
  \end{list}
\end{definition}

\begin{remark}\label{remark_nonsmoothatoms}
  \begin{list}{}{\labelwidth1.3em\leftmargin2em}
    \item[{\upshape (a)\hfill}] The condition \eqref{defcenter} gives us more freedom in choosing the center $x^{\nu,m}$ of each cube $Q_{\nu,m}$. Instead of setting $x^{\nu,m} = 2^{- \nu}m$ as usual, now we can shift the cube $Q_{\nu,m}$ around the point $2^{- \nu}m$ in a range of, at most, $b \, 2^{-\nu}$. This will be useful later on to define cubes on domains.
    \item[{\upshape (b)\hfill}] The number $d$ has the above meaning, see \eqref{d-param}, and it is assumed to be fixed throughout this paper.
    \item[{\upshape (c)\hfill}] As in the smooth case, if $L=0$, then condition (iii) can be ignored since it follows from conditions (i) and (ii) with $K=0$. If $K=0$, then by Definition \ref{holder-space} we only require $a$ to be suitable bounded.
    \item[{\upshape (d)\hfill}] The modification of condition (ii) here was suggested in \cite{TW96} (with some minor adjustments) and it was motivated by the use of the Whitney's extension method to extend atoms from $\mathcal{C}^K(\overline{\Omega})$ to $\mathcal{C}^K(\rn)$, with $0\leq K \notin \nat$. One can see that the usual formulation
        \beq
           \|a(2^{-\nu} \cdot) \mid C^K (\rn) \| \leq c\  \nonumber
        \eeq
        follows from condition (ii) if $K$ is a natural number, since $C^K (\rn) \hookrightarrow \mathcal{C}^K(\rn)$.
    \item[{\upshape (e)\hfill}] Regarding condition (iii), the modification here was suggested by Skrzypczak in \cite{Skr98} for natural numbers $L+1$ (replacing $\mathcal{C}^L(\rn)$ by $C^L(\rn)$). Here, as in \cite{Scharf13}, we extended this definition to general positive numbers $L$.
  \end{list}
\end{remark}

\begin{definition}\label{sequence_spaces}
Let $\bm{w}=(w_\nu)_{\nu\in\nat_0}\in \mathcal{W}^{\alpha}_{\alpha_1,\alpha_2}(\rn)$ and $p,q \in \mathcal{P}(\rn)$.

\begin{list}{}{\labelwidth1.3em\leftmargin2em}
  \item[{\upshape (i)\hfill}] The sequence space $b_{p(\cdot),\q}^{\bm{w}}(\rn)$ consists of those complex-valued sequences $\lambda=(\lambda_{\nu,m})_{\nu\in\nat_0,m\in\zn}$ such that
    \begin{align}
      \|\lambda \mid b_{p(\cdot),\q}^{\bm{w}}(\rn)\|&:= \Big\Vert \Big(\sum_{m\in \zn}|\lambda_{\nu,m}|\,w_{\nu}(2^{-\nu}m)\,\chi_{\nu,m} \Big)_{\nu\in\nat_0}\mid \ell_{\q}(L_{p(\cdot)}(\rn))\Big\Vert \notag
    \end{align}
  is finite.

  \item[{\upshape (ii)\hfill}] If $p^{+}, q^+< \infty$, then the sequence space $f_{p(\cdot),\q}^{\bm{w}}(\rn)$ consists of those complex-valued sequences $\lambda=(\lambda_{\nu,m})_{\nu\in\nat_0,m\in\zn}$ such that
    \begin{align}
      \|\lambda \mid f_{p(\cdot),\q}^{\bm{w}}(\rn)\|&:= \Big\Vert \Big(\sum_{m\in \zn}|\lambda_{\nu,m}|\,w_{\nu}(2^{-\nu}m)\,\chi_{\nu,m} \Big)_{\nu\in\nat_0}\mid L_{p(\cdot)}(\ell_{q(\cdot)} (\rn))\Big\Vert\notag
    \end{align}
    is finite.

\end{list}
\end{definition}

The next result can be found in \cite{GK15} and it states the possibility of decomposing a function $f$ of $B_{\p, \q}^{\bm{w}}(\rn)$ or $F_{\p, \q}^{\bm{w}}(\rn)$ as a linear combination of non-smooth $[K,L]$-atoms according to Definition \ref{defatoms}.

\begin{theorem} \label{atomic_decomp} Let $\bm{w}=(w_{\nu})_{\nu\in\nat_0}\in \mathcal{W}^{\alpha}_{\alpha_1,\alpha_2}(\rn)$ and $p,q \in \mathcal{P}^{log}(\rn)$.
  \begin{list}{}{\labelwidth1.3em\leftmargin2em}
    \item[{\upshape (i)\hfill}] Let $K, L\geq 0$ with $K> \alpha_2$ and $L>\sigma_p - \alpha_1+c_{log}(1/q)$. Then $f \in \mathcal{S}'(\rn)$ belongs to $B_{\p, \q}^{\bm{w}}(\rn)$ if and only if it can be represented as
      \beq
	f=\sum_{\nu=0}^{\infty} \sum_{m\in\zn} \lambda_{\nu,m}\,a_{\nu,m}, \quad\mbox{convergence being in } {\mathcal S}'(\rn) ,
      \eeq
     for $(a_{\nu,m})_{\nu \in\nat_0,m\in\zn}$ non-smooth $[K,L]$-atoms according to Definition \ref{defatoms} and $\lambda\in b_{p(\cdot),\q}^{\bm{w}}(\rn)$. Moreover,
     \beq
	\|f\mid B_{p(\cdot),\q}^{\bm{w}}(\rn)\| \sim \inf \|\lambda \mid b_{p(\cdot),\q}^{\bm{w}}(\rn)\|, \nonumber
     \eeq
     where the infimum is taken over all possible representations of $f$.
     \item[{\upshape (ii)\hfill}] Let $K, L\geq 0$ with $K> \alpha_2$ and $L>\sigma_{p,q} - \alpha_1$. If $p^+,q^+<\infty$, then $f \in \mathcal{S}'(\rn)$ belongs to $F_{\p, \q}^{\bm{w}}(\rn)$ if and only if it can be represented as
      \beq
	f=\sum_{\nu=0}^{\infty} \sum_{m\in\zn} \lambda_{\nu,m}\,a_{\nu,m}, \quad\mbox{convergence being in } {\mathcal S}'(\rn) , \nonumber
      \eeq
     for $(a_{\nu,m})_{\nu \in\nat_0,m\in\zn}$ non-smooth $[K,L]$-atoms according to Definition \ref{defatoms} and $\lambda\in f_{p(\cdot),\q}^{\bm{w}}(\rn)$. Moreover,
     \beq
	\|f\mid F_{p(\cdot),\q}^{\bm{w}}(\rn)\| \sim \inf \|\lambda \mid f_{p(\cdot),\q}^{\bm{w}}(\rn)\|,\nonumber
     \eeq
     where the infimum is taken over all possible representations of $f$.
  \end{list}
\end{theorem}

\begin{remark} This result was proved using the local means characterization of these spaces. The key of the proof is the fact that the local means can also be understood as non-smooth atoms, which allows one to consider estimates of type (ii) and (iii) in Definition \ref{defatoms} in both functions.
\end{remark}

\section{Function spaces on domains}
\subsection{Definitions}

\indent

An open connected set $\Omega$ in $\rn$ is called a domain. As usual $D'(\Omega)$ stands for all complex distributions on the domain $\Omega$ in $\rn$. The restriction of $g \in S'(\rn)$ to $\Omega$ is denoted by $g|_{\Omega}$ and is considered as an element of $D'(\Omega)$.

\begin{definition} \label{defspaces_domains}
Let $\Omega$ be a domain in $\rn$. Let $\bm{w}=(w_{\nu})_{\nu\in\nat_0}\in\mathcal{W}^{\alpha}_{\alpha_1,\alpha_2}(\rn)$ and $p,q \in \mathcal{P}^{\log}(\rn)$.
\begin{list}{}{\labelwidth1.3em\leftmargin2em}
    \item[{\upshape (i)\hfill}] The space $B_{\p,\q}^{\bm{w}}(\Omega)$ is the restriction of $B_{\p,\q}^{\bm{w}}(\rn)$ to $\Omega$, quasi-normed by
        \beq
            \| f \mid B_{\p,\q}^{\bm{w}}(\Omega)\| = \inf \| g \mid B_{p(\cdot),\q}^{\bm{w}}(\rn)\|
        \eeq
        where the infimum is taken over all $g \in B_{p(\cdot),\q}^{\bm{}}(\rn)$ with $g|_{\Omega}=f$.
     \item[{\upshape (ii)\hfill}] If $p^{+}, q^+< \infty$, then the space $F_{\p,\q}^{\bm{w}}(\Omega)$ is the restriction of $F_{\p,\q)}^{\bm{w}}(\rn)$ to $\Omega$, quasi-normed by
        \beq
            \| f \mid F_{\p,\q}^{\bm{w}}(\Omega)\| = \inf \| g \mid F_{p(\cdot),\q}^{\bm{w}}(\rn)\|
        \eeq
        where the infimum is taken over all $g \in F_{p(\cdot),\q}^{\bm{w}}(\rn)$ with $g|_{\Omega}=f$.\\
\end{list}
\end{definition}

\begin{remark}\label{rem:extension} The definition requires that $\p, \q$ and $\bm{w}$ are defined on all of $\rn$. Working with the restrictions above we only need to consider the values of $\p, \q$ and $\bm{w}$ on $\Omega$. On the other hand, one can directly start with functions $\p, \q$ and $\bm{w}$ which are only defined on $\Omega$ and extend them to whole of $\rn$. Such extensions are usually not unique, but for example the extension of functions from $C^{\log}(\Omega)$ to $C^{\log}(\rn)$ in \cite[Proposition 4.1.7]{DHHR} preserves fundamental properties. 

To the best of the authors knowledge, it is unknown if the restricted spaces of Definition \ref{defspaces_domains} are independent of the extension of the parameters from $\Omega$ to $\rn$. It remains, therefore, an open problem which may be considered in future work.
\end{remark}

\begin{definition}
Let $MR(n)$ (minimally regular) be the collection of all bounded domains $\Omega$ in $\rn$ with
\beq
    \Omega = int (\overline{\Omega}),
\eeq
that means, $\Omega$ coincides with the interior of its closure $\overline{\Omega}$.
\end{definition}

\begin{remark} For more details regarding these domains, we refer section 3.1 in \cite{TW96}.
\end{remark}

\subsection{Regular domains}

\indent

Our aim is to characterize the spaces $B_{\p,\q}^{\bm{w}}(\Omega)$ and $F_{\p,\q}^{\bm{w}}(\Omega)$ intrinsically using non-smooth atoms. To this end we resort to the types of domains already used by Triebel and Winkelvoß in \cite{TW96}, where intrinsic atomic characterizations of the classical spaces $B^s_{p,q}(\Omega)$ and $F^s_{p,q}(\Omega)$ were found. We describe now these domains, which are naturally connected with our task. Let $\partial \Omega = \overline{\Omega} \setminus \Omega$ denote the boundary of $\Omega$.
\begin{definition}
    \begin{list}{}{\labelwidth1.3em\leftmargin2em}
    \item[{\upshape (i)\hfill}] Let $IR(n)$ (interior regular) be the collection of all domains $\Omega \in MR(n)$ for which one finds a positive number $c$ such that for any cube $Q$ centered at $\partial \Omega$ with side-length less than or equal 1,
        \beq \label{interior_reg} |Q \cap \Omega| \geq c \, |Q|.
        \eeq
    \item[{\upshape (ii)\hfill}] Let $ER(n)$ (exterior regular) be the collection of all domains $\Omega \in MR(n)$ for which one finds a positive number $c$ such that any cube $Q$ centered at $\partial \Omega$ with side-length $l$ less than or equal 1, there exists a subcube $Q^e$ with side-length $c l$ and
        \beq \label{exterior_reg}
            Q^e \subset Q \cap (\rn \setminus \overline{\Omega}).
        \eeq
    \item[{\upshape (iii)\hfill}]  Let
        \beq   \label{regular}
            R(n) = IR(n) \cap ER(n)
        \eeq
        be the collection of all domains $\Omega \in MR(n)$ which are both interior and exterior regular.
    \end{list}
\end{definition}

\begin{remark}
\begin{list}{}{\labelwidth1.3em\leftmargin2em}
    \item[{\upshape (a)\hfill}] Analogously to $ER(n)$, let $\Omega \in MR(n)$ be a domain for which one finds a positive number $c$ such that for any cube $Q$ centered at $\partial \Omega$ with side-length $l$ less than or equal 1, there exists a subcube $Q^i$ with side length $cl$ and
            \beq \label{interior_reg1}
                Q^i \subset Q \cap \Omega.
            \eeq
        Then we have $\Omega \in IR(n)$. However, although this condition is quite natural in order to have $\Omega \in IR(n)$, there are domains $\Omega \in IR(n)$ for which \eqref{interior_reg1} is not true. If one takes out of a square in $\mathbb{R}^2$ infinitely many smaller squares such that one obtains a carpet-like domain, then it might happen that \eqref{interior_reg1} is violated but not \eqref{interior_reg}.
    \item[{\upshape (b)\hfill}] Regarding specific (non-smooth) domains connected with these definitions, we mention that if $\Omega \in MR(n)$ is a so called $(\epsilon,\delta)$-domain, then it belongs to $IR(n)$. In particular, any connected bounded Lipschitz domain is an interior regular domain. For more details, we refer \cite{TW96}. Moreover, in \cite{HKT08} the authors considered domains satisfying the \textit{measure density condition}, which actually coincide with our definition of interior regular domain.
    \item[{\upshape (c)\hfill}] Similar but not identical with this class of domains is the class of \textit{thick domains}, $E$-thick and $I$-thick. For more details, see \cite{Tr08}.
\end{list}
\end{remark}

\subsection{Atoms on domains}

\indent

We always assume $\Omega \in MR(n)$.

\begin{definition} Let $s> 0$ and $\Omega \in MR(n)$. Then $\mathcal{C}^{s}(\overline{\Omega})$ consists of all complex-valued continuous functions $f$ on $\overline{\Omega}$ with the following two properties:
\begin{list}{}{\labelwidth1.3em\leftmargin2em}
    \item[{\upshape (i)\hfill}] $f$ has classical derivatives $D^{\alpha}f$ in $\Omega$ for $|\alpha| \leq \lfloor s\rfloor^-$ and there exist continuous functions $f_{\alpha}$ on $\overline{\Omega}$ which coincide with $D^{\alpha}f$ on $\Omega$,
    \item[{\upshape (ii) \hfill}]
        \beq
             \| f \mid \mathcal{C}^s(\overline{\Omega}) \|:= \sum_{|\alpha|\leq \lfloor s\rfloor^-} \sup_{x \in \overline{\Omega}} |D^{\alpha}f(x)| + \sum_{|\alpha| = \lfloor s\rfloor^-} \sup_{x,y \in \overline{\Omega}, x\neq y} \frac{|f(x)-f(y)|}{|x-y|^{\{s\}^+}} < \infty. \nonumber
        \eeq
\end{list}
\end{definition}

In order to introduce atoms on domains $\Omega \in MR(n)$ we again rely on the cubes $Q_{\nu, m}$. We may assume in addition that the centers $x^{\nu,m}$ of the cubes $Q_{\nu,m}$ with $d\, Q_{\nu,m} \cap \partial\Omega \neq \emptyset$ are located at $\partial \Omega$. In this sense we call $Q_{\nu, m}$
\beq \label{interior_cube}
    \mbox{an interior cube if } d\, Q_{\nu, m} \subset \Omega, \qquad \nu \in \nat_0, m \in \zn,
\eeq
and
\beq\label{boundary_cube}
    \mbox{a boundary cube if } x^{\nu,m} \in  \partial\Omega, \qquad \nu \in \nat_0, m \in \zn.
\eeq
Other cubes are not of interest for us. Let, for brevity,
\beq
    \Omega^{\nu} = \{ x \in \rn : 2^{-\nu}x \in \Omega\}, \qquad \nu \in \nat_0.
\eeq

\begin{definition} \label{def_non-smooth_atoms} Let $\Omega \in MR(n)$, $d>1$ and $c>0$.

    \begin{list}{}{\labelwidth1.3em\leftmargin2em}
        \item[{\upshape (a)\hfill}] Let $K, L\geq 0$.  Then $a(x)$ is called a non-smooth interior $[K, L]$-atom in $\Omega$, for all $\nu \in \nat_0$ and $m \in \zn$, if
            \begin{list}{}{\labelwidth1.3em\leftmargin2em}
                \item[{\upshape (i)\hfill}] $\supp a \subset d\,Q_{\nu,m}, \quad$ for some interior cube $Q_{\nu, m}$,
                \item[{\upshape (ii)\hfill}] $ \|a(2^{-\nu} \cdot) \mid \mathcal{C}^K (\overline{\Omega^{\nu}}) \| \leq c\,$,
                \item[{\upshape (iii)\hfill}] and for every $\psi \in \mathcal{C}^L(\rn)$ it holds
                    \beq \label{moment_cond}
                        \Big| \ds \int_{d Q_{\nu, m}} \psi(x) a(x) dx\Big| \leq c \, 2^{-\nu(L + n)}\| \psi \mid \mathcal{C}^L(\rn)\|.
                    \eeq
            \end{list}

        \item[{\upshape (b)\hfill}] Let $K \geq 0$. Then $a(x)$ is called a non-smooth boundary $[K, 0]$-atom in $\Omega$, for all $\nu \in \nat_0$ and $m \in \zn$, if
            \begin{list}{}{\labelwidth1.3em\leftmargin2em}
                \item[{\upshape (i)\hfill}] $\supp a \subset \overline{\Omega} \cap d\,Q_{\nu,m}, \quad$  for some boundary cube $Q_{\nu, m}$,
                \item[{\upshape (ii)\hfill}] $ \|a(2^{-\nu} \cdot) \mid \mathcal{C}^K (\overline{\Omega^{\nu}}) \| \leq c\,$.
            \end{list}

    \end{list}
\end{definition}
\begin{remark}
The above part (a) is the natural counterpart of Definition \ref{defatoms}. As for part (b) no conditions of type \eqref{moment_cond} are required.
\end{remark}

\subsection{Atomic domains}

\indent

We start by introducing the counterparts of the sequence spaces $b_{\p,\q}^{\bm{w}}$ and $f_{\p,\q}^{\bm{w}}$ from Definition \ref{sequence_spaces}. Let $\Omega \in MR(n)$ and let $Q_{\nu,m}$ be the dyadic cubes defined above, where we are only interested in interior and boundary cubes described in \eqref{interior_cube} and \eqref{boundary_cube}, respectively. Let
\beq \label{sequences_collection}
    \lambda = \{ \lambda_{\nu,m} : \lambda_{\nu,m} \in \mathbb{C}, \nu \in \nat_0, m \in \zn, Q_{\nu,m} \mbox{ interior or boundary cube}\}.
\eeq
Furthermore, $\displaystyle \sum_{m \in \zn} {\!}^{\nu, \Omega}$ means that for fixed $\nu \in \nat_0$ the sum is taken over those $m \in \zn$ for which $Q_{\nu,m}$ is an interior or boundary cube.

\begin{definition}\label{sequence_spaces_domains}
Let $\Omega \in MR(n)$, $\bm{w}=(w_\nu)_{\nu\in\nat_0}\in \mathcal{W}^{\alpha}_{\alpha_1,\alpha_2}(\rn)$ and $p,q \in \mathcal{P}(\rn)$.

\begin{list}{}{\labelwidth1.3em\leftmargin2em}
  \item[{\upshape (i)\hfill}] The sequence space $b_{\p,\q}^{\bm{w}}(\Omega)$ consists of those complex-valued sequences $\lambda=(\lambda_{\nu,m})_{\nu\in\nat_0,m\in\zn}$ given by \eqref{sequences_collection} such that
    \begin{align}
      \|\lambda \mid b_{\p,\q}^{\bm{w}}(\Omega)\|&:= \Big\Vert \Big(\sum_{m\in \zn}{\!}^{\nu, \Omega}\,|\lambda_{\nu,m}|\,w_{\nu}(2^{-\nu}m)\,\chi_{\nu,m} \Big)_{\nu\in\nat_0}\mid \ell_{\q}(L_{\p}(\Omega))\Big\Vert \notag
    \end{align}
  is finite.

  \item[{\upshape (ii)\hfill}] If $p^{+}, q^+< \infty$, then the sequence space $f_{\p,\q}^{\bm{w}}(\Omega)$ consists of those complex-valued sequences $\lambda=(\lambda_{\nu,m})_{\nu\in\nat_0,m\in\zn}$ given by \eqref{sequences_collection} such that
    \begin{align}
      \|\lambda \mid f_{\p,\q}^{\bm{w}}(\Omega)\|&:= \Big\Vert \Big(\sum_{m\in \zn}{\!}^{\nu, \Omega}\,|\lambda_{\nu,m}|\,w_{\nu}(2^{-\nu}m)\,\chi_{\nu,m} \Big)_{\nu\in\nat_0}\mid L_{\p}(\ell_{\q} (\Omega))\Big\Vert\notag
    \end{align}
    is finite.

\end{list}
\end{definition}

Next we are interested in the counterpart of Theorem \ref{atomic_decomp}. Regarding the question whether the corresponding series converges, we convert it now in a definition of domains having this property. The conditions which appear in this definition are the natural restrictions on the parameters $K$ and $L$ that we already have in Theorem \ref{atomic_decomp}. This justifies to take over this knowledge to the situation we consider now.

\begin{definition}\label{atomic_domains}
    Let $\bm{w}=(w_\nu)_{\nu\in\nat_0}\in \mathcal{W}^{\alpha}_{\alpha_1,\alpha_2}(\rn)$ and $p,q \in \mathcal{P}^{log}(\rn)$.
  \begin{list}{}{\labelwidth1.3em\leftmargin2em}
    \item[{\upshape (i)\hfill}] Let $K, L\geq 0$ with $K> \alpha_2$ and $L>\sigma_{p} - \alpha_1+c_{log}(1/q)$. Then $ Atom(B_{\p, \q}^{\bm{w}})^n$ (atomic $B_{\p, \q}^{\bm{w}}$-domain) denotes the collection of all domains $\Omega \in MR(n)$ such that for all choices of $K$ and $L$
      \beq \label{series0}
	\sum_{\nu=0}^{\infty} \sum_{m\in\zn}{\!}^{\nu, \Omega}\, \lambda_{\nu,m}\,a_{\nu,m}, \quad \lambda \in b_{\p,\q}^{\bm{w}}(\Omega),
      \eeq
     converges in $D'(\Omega)$ to an element of $B_{\p, \q}^{\bm{w}}(\Omega)$, where $(a_{\nu,m})_{\nu \in\nat_0,m\in\zn}$ are non-smooth interior $[K,L]$-atoms, or non-smooth boundary $[K,0]$-atoms according to Definition \ref{def_non-smooth_atoms}.

     \item[{\upshape (ii)\hfill}] Let $K, L\geq 0$ with $K> \alpha_2$ and $L>\sigma_{p,q} - \alpha_1$. If $p^+,q^+<\infty$, then $ Atom(F_{\p, \q}^{\bm{w}})^n$ (atomic $F_{\p, \q}^{\bm{w}}$-domain) denotes the collection of all domains $\Omega \in MR(n)$ such that for all choices of $K$ and $L$
      \beq
	   \sum_{\nu=0}^{\infty} \sum_{m\in\zn}{\!}^{\nu, \Omega}\, \lambda_{\nu,m}\,a_{\nu,m}, \quad \lambda \in f_{\p, \q}^{\bm{w}}(\Omega),
      \eeq
     converges in $D'(\Omega)$ to an element of $F_{\p, \q}^{\bm{w}}(\Omega)$, where $(a_{\nu,m})_{\nu \in\nat_0,m\in\zn}$ are non-smooth interior $[K,L]$-atoms, or non-smooth boundary $[K,0]$-atoms according to Definition \ref{def_non-smooth_atoms}.
  \end{list}
\end{definition}

\subsection{Atomic characterizations}

\indent

Before stating the main result, we present two lemmas which will be useful in the sequel. Similarly to \cite[Lemma 1.19]{Mou01}, the first lemma states that, for fixed $\nu \in \nat_0$, each $x\in \rn$ belongs to a finite number of cubes $d Q_{\nu,m}$. We make use of this result and the additional Lemmas \ref{conv_1} and \ref{conv_2} to prove the second helpful result, stated in Lemma \ref{Q->E}. It shows that we get an equivalent norm for the sequence spaces $b_{\p, \q}^{\bm{w}}(\rn)$ and $f_{\p, \q}^{\bm{w}}(\rn)$ when we shift a bit around the cubes $Q_{\nu,m}$. For cubes with center $x^{\nu,m} = 2^{-\nu}m$, this result was already proven in \cite{Mou13} and \cite{GMN14}. Here we adapt the proof to our needs.

\begin{lemma}\label{finite_cubes}
Fix $\nu \in \nat_0$ and let $b$ and $d$ be as before. Then any $x \in \rn$ belongs to at most $N$ cubes $d Q_{\nu,m}, m\in \zn$, where $N$ is independent of $\nu$ and $m$ (it only depends on $b,d$ and on the dimension $n$).
\end{lemma}
\bpr
    Let $x\in \rn$. By \eqref{d-param} there surely exists $m \in \zn$ such that $x \in d Q_{\nu,m}$, which means that
    \beq
	|x_i-x_i^{\nu,m}| \leq d 2^{-\nu -1}, \quad i=1, ... , n, \nonumber
    \eeq
    or equivalently, by \eqref{defcenter}, that
    \beq
	|2^{\nu}x_i - m_i| \leq b + \frac{d}{2}, \quad i=1, ... , n. \nonumber
    \eeq

    Assume that $x\in dQ_{\nu, m'}$ for some $m'\in \zn$ with $m'\neq m$. Similarly as before, we get that
    \beq
	|2^{\nu}x_i - m_i'| \leq b + \frac{d}{2}, \quad i=1, ... , n. \nonumber
    \eeq
    This gives
    \beq
	|m_i-m_i'| \leq 2b+d, \quad i=1, ..., n, \nonumber
    \eeq
    which means that $m'$ belongs to the cube centered at $m$ and with side length $2(2b+d)$. The number of such $m' \in \zn$ is $N= 2^n[2b+d]^n$, where $[a]$ denotes the integer part of $a$.
\epr

The next two results state convolution inequalities for $B_{p(\cdot),\q}^{\bm{w}}(\rn)$ and $F_{p(\cdot),\q}^{\bm{w}}(\rn)$ and can be found in \cite{KemV12} and \cite{DHR09}, respectively. We note that a slightly different version of Lemma \ref{conv_1} was firstly proved in \cite[Lemma 4.7]{AH10}. We introduce the functions
$$ \eta_{\nu, R}(x) = \frac{2^{n \nu }}{(1+2^{\nu }|x|)^R},$$
for $\nu \in \nat_0$ and $R>0$.

\begin{lemma}[{{\!\!\cite[Lemma 10]{KemV12}}}]\label{conv_1} Let $p, q \in \mathcal{P}(\rn)$ with $\p \geq 1$. For all $R>n +c_{log}(1/q)$, there exists a constant $c>0$ such that for all sequences $(f_{\nu})_{\nu \in \nat_0} \in \ell_{\q}(L_{\p}(\rn))$ it holds
\beq
  \| \eta_{\nu, R} \ast f_{\nu} \mid \ell_{\q}(L_{\p}(\rn)) \| \leq c \, \|f_{\nu} \mid \ell_{\q}(L_{\p}(\rn)) \|. \nonumber
\eeq
\end{lemma}

\begin{lemma}[{{\!\!\cite[Theorem 3.2]{DHR09}}}]\label{conv_2} Let $p, q \in C^{log}(\rn)$ with $1<p^-\leq p^+ <\infty$ and $1<q^-\leq q^+ <\infty$. Then the inequality
\beq
  \| \eta_{\nu, R} \ast f_{\nu} \mid L_{\p}(\ell_{\q}(\rn)) \| \leq c\, \|f_{\nu} \mid L_{\p}(\ell_{\q}(\rn)) \| \nonumber
\eeq
holds for every sequence $(f_{\nu})_{\nu \in \nat_0}$ of $L_1^{loc}(\rn)$ functions and constant $R>n$.
\end{lemma}

\begin{lemma}\label{Q->E}
Let $\bm{w}=(w_\nu)_{\nu\in\nat_0}\in\mathcal{W}^{\alpha}_{\alpha_1,\alpha_2}(\rn)$, $p, q\in \mathcal{P}^{\log}(\rn)$.  Let $d, \varepsilon>0$ and let $\{E_{\nu,m}\}_{\nu\in\nat_0,m\in\zn}$ be a collection of measurable sets with $E_{\nu,m}\subset dQ_{\nu,m}$ and $|E_{\nu,m}|\geq \varepsilon |Q_{\nu,m}|$, for all $\nu\in\nat_0$ and $m\in \zn$.
\begin{list}{}{\labelwidth1.3em\leftmargin2em}
    \item[{\upshape (i)\hfill}] Then
        $$\|\lambda \mid b^{\bm{w}}_{\p,\q}(\rn) \| \sim \Big\Vert \Big(\displaystyle \sum_{m\in\mathbb{Z}^n} |\lambda_{\nu,m}| w_{\nu}(2^{-\nu}m) \chi_{E_{\nu,m}}\Big)_{\nu \in \nat_0} \mid \ell_{\q}(L_{\p}(\rn))\Big\Vert.$$
    \item[{\upshape (ii)\hfill}] If $p^{+}, q^+<\infty$, then
        $$\|\lambda \mid f^{\bm{w}}_{\p,\q}(\rn) \| \sim \Big\Vert \Big(\displaystyle \sum_{m\in\mathbb{Z}^n} |\lambda_{\nu,m}| w_{\nu}(2^{-\nu}m) \chi_{E_{\nu,m}}\Big)_{\nu \in \nat_0} \mid L_{\p}(\ell_{\q}(\rn))\Big\Vert.$$
\end{list}
\end{lemma}
\bpr
    We will present the proof of (i) since the other case follows similarly. Starting with the inequality ``$\leq$'', let $0<r<\min(1,p^-)$. We express the norm as
\begin{align*}
\|\lambda \mid b_{p(\cdot),\q}^{\bm{w}}(\rn)\|
&=
\,\Big\Vert \Big(\sum_{m\in \zn}|\lambda_{\nu,m}|\,w_{\nu}(2^{-\nu}m)\,\chi_{\nu,m} \Big)_{j\in\nat_0}\mid \ell_{\q}(L_{\p}(\rn))\Big\Vert\\
&=
\,\Big\Vert\Big(\sum_{m\in \zn}|\lambda_{\nu,m}| \, w_{\nu}(2^{-\nu}m)\,\chi_{\nu,m}\Big)^r_{\nu\in\nat_0}\mid \ell_{\frac{\q}{r}}(L_{\frac{\p}{r}}(\rn))\Big\Vert^{\frac 1r}\\
&\leq\,\Big\Vert\Big(\sum_{m\in \zn}|\lambda_{\nu,m}|^r \, w_{\nu}^r(2^{-\nu}m)\,\chi_{\nu,m}\Big)_{\nu\in\nat_0}\mid \ell_{\frac{\q}{r}}(L_{\frac{\p}{r}}(\rn))\Big\Vert^{\frac 1r}
\end{align*}
where the last step is true by Lemma \ref{finite_cubes} and the fact that $r <1$. Now, for each $R>0$, we use the estimate $\chi_{\nu,m}\leq c\eta_{\nu,R}\ast \chi_{E_{\nu,m}}$ for all $\nu\in\nat_0$ and $m\in\zn$. Choosing $R>n+c_{log}(1/q)$, we use Lemma \ref{conv_1} to derive the following:
\begin{align*}
\|\lambda \mid b_{\p,\q}^{\bm{w}}(\rn)\|
&\leq c
\,\Big\Vert\Big(\sum_{m\in \zn}|\lambda_{\nu,m}|^r\,w_{\nu}^r(2^{-\nu}m)\, \eta_{\nu,R} \ast \chi_{E_{\nu,m}}\Big)_{\nu\in\nat_0}\mid \ell_{\frac{\q}{r}}(L_{\frac{\p}{r}}(\rn))\Big\Vert^{\frac 1r}\\
&  = c
\,\Big\Vert\Big( \eta_{\nu,R}\ast \Big(\sum_{m\in \zn}|\lambda_{\nu,m}|^r\,w_{\nu}^r(2^{-\nu}m)\, \chi_{E_{\nu,m}}\Big)\Big)_{\nu\in\nat_0}\mid \ell_{\frac{\q}{r}}(L_{\frac{\p}{r}}(\rn))\Big\Vert^{\frac 1r}\\
&\leq c'
\,\Big\Vert\Big(\sum_{m\in \zn}|\lambda_{\nu,m}|^r\,w_{\nu}^r(2^{-\nu}m)\, \chi_{E_{\nu,m}}\Big)_{\nu\in\nat_0}\mid \ell_{\frac{\q}{r}}(L_{\frac{\p}{r}}(\rn))\Big\Vert^{\frac 1r}\\
& \sim
\Big\Vert\Big(\sum_{m\in \zn}|\lambda_{\nu,m}|\,w_{\nu}(2^{-\nu}m)\, \chi_{E_{\nu,m}}\Big)_{\nu  \in \mathbb{N}_0}\mid \ell_{\q}(L_{p(\cdot)}(\rn))\Big\Vert.
\end{align*}
The other direction follows by the same arguments since, for fixed $R>0$,  $\chi_{E_{\nu,m}}\leq c\eta_{\nu,R}\ast \chi_{\nu,m}$ for all $\nu\in\nat_0$ and $m\in\zn$. We use in this case Lemma \ref{conv_2}.
\epr

If one wishes to extend boundary atoms in the sense of Definition \ref{def_non-smooth_atoms} (b) from $\Omega$ to $\rn$, then the conditions of type (iii) in Definition \ref{defatoms} cause some trouble. So we start by avoiding this problem stating a result where $L=0$ can be considered. For the interior atoms in the sense of Definition \ref{def_non-smooth_atoms}, we keep all the conditions.
\begin{theorem}
    \label{atomic_decomp_domains} Let $\bm{w}=(w_{\nu})_{\nu\in\nat_0}\in \mathcal{W}^{\alpha}_{\alpha_1,\alpha_2}(\rn)$ and $p,q \in \mathcal{P}^{log}(\rn)$.
  \begin{list}{}{\labelwidth1.3em\leftmargin2em}
    \item[{\upshape (i)\hfill}] Let $\alpha_1 >\sigma_{p} + c_{log}(1/q)$. Then
        \beq \label{atomicdomain_B} Atom(B_{\p,\q}^{\bm{w}})^n \supset IR(n).
        \eeq

    Let $K, L\geq 0$ with $K> \alpha_2$. Then $f \in \mathcal{D}'(\Omega)$ belongs to $B_{\p, \q}^{\bm{w}}(\Omega)$ if and only if it can be represented as
      \beq \label{series1}
	f=\sum_{\nu=0}^{\infty} \sum_{m\in\zn}{\!}^{\nu, \Omega}\, \lambda_{\nu,m}\,a_{\nu,m}, \quad\mbox{convergence being in } {\mathcal D}'(\Omega) ,
      \eeq
     in the sense of Definition \ref{atomic_domains}, where $(a_{\nu,m})_{\nu \in\nat_0,m\in\zn}$ are non-smooth interior $[K,L]$-atoms, or non-smooth boundary $[K,0]$-atoms according to Definition \ref{def_non-smooth_atoms} and $\lambda\in b_{\p, \q}^{\bm{w}}(\Omega)$. Moreover,
     \beq
	\|f\mid B_{\p, \q}^{\bm{w}}(\Omega)\| \sim \inf \|\lambda \mid b_{\p, \q}^{\bm{w}}(\Omega)\|, \nonumber
     \eeq
     where the infimum is taken over all possible representations of $f$.
     \item[{\upshape (ii)\hfill}] Let $\alpha_1 >\sigma_{p,q}$ and $p^+,q^+<\infty$. Then
        \beq Atom(F_{\p,\q}^{\bm{w}})^n \supset IR(n).
        \eeq

      Let $K, L\geq 0$ with $K> \alpha_2$. Then $f \in \mathcal{D}'(\Omega)$ belongs to $F_{\p,\q}^{\bm{w}}(\Omega)$ if and only if it can be represented as
      \beq \label{series2}
	f=\sum_{\nu=0}^{\infty} \sum_{m\in\zn}{\!}^{\nu, \Omega}\, \lambda_{\nu,m}\,a_{\nu,m}, \quad\mbox{convergence being in } {\mathcal D}'(\Omega) ,
      \eeq
     in the sense of Definition \ref{atomic_domains}, where $(a_{\nu,m})_{\nu \in\nat_0,m\in\zn}$ are non-smooth interior $[K,L]$-atoms, or non-smooth boundary $[K,0]$-atoms according to Definition \ref{def_non-smooth_atoms} and $\lambda\in f_{\p,\q}^{\bm{w}}(\Omega)$. Moreover,
     \beq
	\|f\mid F_{\p,\q}^{\bm{w}}(\Omega)\| \sim \inf \|\lambda \mid f_{\p,\q}^{\bm{w}}(\Omega)\|,\nonumber
     \eeq
     where the infimum is taken over all possible representations of $f$.
  \end{list}
\end{theorem}
\bpr
\underline{Step 1.} We only prove part (i), since part (ii) follows similarly. Let $\Omega \in IR(n)$. First we want to show that the series \eqref{series1} converges in $\mathcal{D}'(\Omega)$. For this purpose we extend each atom $a_{\nu,m}(x)$ individually from $\Omega$ to $\rn$. We follow the proof of \cite[Theorem 3.5]{TW96}, which relies on Whitney's extension method, according to \cite{Stein}, pp. 170-180. Note that only the boundary atoms in the sense of Definition \ref{def_non-smooth_atoms} (b) are of interest. Let $\mathcal{E}_k$ be the linear extension operator constructed in \cite{Stein}, p. 177, formula (18). By Theorem 4 on the same page, $\mathcal{E}_{ [ \sigma ]}$ generates a linear extension operator
\beq
    \mathcal{E}_{ [ \sigma ]}(\overline{\Omega}) : \mathcal{C}^{\sigma}(\overline{\Omega}) \rightarrow \mathcal{C}^{\sigma}(\rn), \nonumber
\eeq
with a bound being independent of $\overline{\Omega}$. Let $D_c$ be the dilation operator on $\rn$,
\beq
    D_c : f(x) \mapsto f(cx), \quad c>0. \nonumber
\eeq
Then, it follows from the explicit construction of $\mathcal{E}_{ [ \sigma ] }$ that
\beq
    \mathcal{E}_{ [ \sigma ]}(\overline{\Omega}) = D_{2^{\nu}} \circ \mathcal{E}_{ [ \sigma ]}(\overline{\Omega^{\nu}}) \circ D_{2^{-\nu}}, \quad \nu \in \nat_0. \nonumber
\eeq
Let $\psi$ be a $C^{\infty}$ cut-off function with
\beq
    \supp \psi \subset 2dQ \quad \mbox{ and } \quad \psi(x)=1 \mbox{ if } x \in dQ, \nonumber
\eeq
where $Q$ is the unit cube centered at the origin and $d$ has the same meaning as in Definition \ref{def_non-smooth_atoms} (b). We apply the operator $\mathcal{E}_{ [ \sigma ]}(\overline{\Omega})$ to the non-smooth boundary $[K,0]$-atom $a_{\nu, m}$ and put
\beq
  b_{\nu,m}(x) = \psi(2^{\nu}(x-x^{\nu,m})) (\mathcal{E}_{ [ \sigma ]}(\overline{\Omega})a_{\nu,m})(x). \nonumber
\eeq
In this way, we get new non-smooth atoms $b_{\nu,m}$ on $\rn$ in the sense of Definition \ref{defatoms}, with $L=0$ and after replacing $d$ by $2d$. Then the counterpart of \eqref{series1} is given by
\beq \label{series_counterpart}
    \sum_{\nu=0}^{\infty} \sum_{m\in\zn} \mu_{\nu,m}\,b_{\nu,m}, \quad x \in \rn, \quad \mu \in b_{p(\cdot),q(\cdot)}^{\bm{w}}(\rn),
\eeq
with $\mu_{\nu,m}= b_{\nu,m}(x)=0$ if $Q_{\nu,m}$ is an "exterior" cube and $\mu_{\nu,m}=\lambda_{\nu,m}$ otherwise. Now Theorem \ref{atomic_decomp} allows one  to conclude that the series \eqref{series_counterpart} converges in $\mathcal{S}'(\rn)$ and its limit belongs to $B_{\p,\q}^{\bm{w}}(\rn)$. Hence, its restriction \eqref{series0} to $\Omega$ converges in $\mathcal{D}'(\Omega)$ and, by Definition \ref{defspaces_domains}, its limit belongs to $B_{\p,\q}^{\bm{w}}(\Omega)$. Thus the proof of \eqref{atomicdomain_B} is complete.

\underline{Step 2.} Let us prove now the remaining part of (i). By the above argument, any $f \in \mathcal{D}'(\Omega)$ given by \eqref{series1} with $\lambda \in b_{\p,\q}^{\bm{w}}(\Omega)$ belongs to $B_{\p,\q}^{\bm{w}}(\Omega)$. Conversely, by Definition \ref{defspaces_domains} and Theorem \ref{atomic_decomp}, any $f \in B_{\p,\q}^{\bm{w}}(\Omega)$ can be represented as in \eqref{series1}. Regarding the equivalence of quasi-norms, it follows from Definition \ref{defspaces_domains} and Theorem \ref{atomic_decomp}, aside the fact that
\beq \label{ineq1}
    \|\mu \mid b_{p(\cdot),q(\cdot)}^{\bm{w}}(\rn) \| \leq c \, \|\lambda \mid b_{\p,\q}^{\bm{w}}(\Omega) \|,
\eeq
with $\lambda$ and $\mu$ as before. We prove now inequality \eqref{ineq1}. Let $Q_{\nu,m}$ be a boundary cube with center $x^{\nu,m} \in \partial \Omega$. Since $\Omega \in IR(n)$, we can find $\epsilon > 0$ such that
$$ |Q_{\nu,m} \cap \Omega| \geq \epsilon|Q_{\nu,m}|.$$
Denote by $E_{\nu,m}:= Q_{\nu,m} \cap \Omega$. Then, by Lemma \ref{Q->E}, we have
\begin{align}
     \|\mu \mid b_{p(\cdot),q(\cdot)}^{\bm{w}}(\rn) \| &= \Big\Vert \Big(\sum_{m\in \zn}|\mu_{\nu,m}|\,w_{\nu}(2^{-\nu}m)\,\chi_{\nu,m} \Big)_{\nu\in\nat_0}\mid \ell_{\q}(L_{p(\cdot)}(\rn))\Big\Vert \nonumber \\
     &= \Big\Vert \Big(\sum_{m\in \zn}{\!}^{\nu, \Omega}\, |\lambda_{\nu,m}|\,w_{\nu}(2^{-\nu}m)\,\chi_{\nu,m} \Big)_{\nu\in\nat_0}\mid \ell_{\q}(L_{p(\cdot)}(\rn))\Big\Vert \nonumber \\
     &\sim \Big\Vert \Big(\sum_{m\in \zn}{\!}^{\nu, \Omega}\, |\lambda_{\nu,m}|\,w_{\nu}(2^{-\nu}m)\,\chi_{E_{\nu,m}} \Big)_{\nu\in\nat_0}\mid \ell_{\q}(L_{p(\cdot)}(\rn))\Big\Vert \nonumber \\
     &\sim \Big\Vert \Big(\sum_{m\in \zn}{\!}^{\nu, \Omega}\, |\lambda_{\nu,m}|\,w_{\nu}(2^{-\nu}m)\,\chi_{\nu,m} \Big)_{\nu\in\nat_0}\mid \ell_{\q}(L_{\p}(\Omega))\Big\Vert \nonumber \\
     &= \|\lambda \mid b_{\p,\q}^{\bm{w}}(\Omega) \|. \nonumber
\end{align}
\epr

\begin{remark}
    \begin{list}{}{\labelwidth1.3em\leftmargin2em}
        \item[{\upshape (a)\hfill}] In the case of constant exponents, \cite[Theorem 3.5]{TW96} shows that one only needs $\Omega \in MR(n)$ for $B^s_{p,q}$, whereas for $F^s_{p,q}$ one requires $\Omega \in IR(n)$. This supports the well-known fact that the spaces of type $B^s_{p,q}$ are structurally simpler than the spaces of type $F^s_{p,q}$.\\
				For the variable exponent spaces which are the subject we consider here this situation changes. To prove inequality \eqref{ineq1}, we also need to assume in the case of  $B_{\p,\q}^{\bm{w}}$ that $\Omega \in IR(n)$. 				
        \item[{\upshape (b)\hfill}] Another significant difference from the classical case is the proof of inequality \eqref{ineq1}. For this purpose, Triebel and Winkelvoß used the Hardy-Littlewood maximal operator. Instead we base our proof on Lemma \ref{Q->E}, which comes out through the convolution inequalities stated in Lemmas \ref{conv_1} and \ref{conv_2}.
    \end{list}
\end{remark}

In Theorem \ref{atomic_decomp_domains} no moment conditions of type (iii) in Definition \ref{def_non-smooth_atoms} (a) are required. In this case Whitney's extension method proved to be an effective tool to extend atoms from $\Omega$ (better $\overline{\Omega}$) to $\rn$. When in case of needing such conditions, one has to complement what was done so far in the proof of Theorem \ref{atomic_decomp_domains} by a special method which creates moment conditions on $\rn \backslash \overline{\Omega}$. For this purpose, we need the additional assumption $\Omega \in ER(n)$, as we will present in the next result.

\begin{theorem}
    \label{atomic_decomp_domains-general} Let $\bm{w}=(w_{\nu})_{\nu\in\nat_0}\in \mathcal{W}^{\alpha}_{\alpha_1,\alpha_2}(\rn)$ and $p,q \in \mathcal{P}^{log}(\rn)$.
  \begin{list}{}{\labelwidth1.3em\leftmargin2em}
    \item[{\upshape (i)\hfill}] Then
        \beq Atom(B_{\p, \q}^{\bm{w}})^n \supset R(n).
        \eeq

    Let $K, L\geq 0$ with $K> \alpha_2$ and $L> \sigma_{p} - \alpha_1 + c_{log}(1/q)$. Then $f \in \mathcal{D}'(\Omega)$ belongs to $B_{\p,\q}^{\bm{w}}(\Omega)$ if and only if it can be represented as
      \beq
	f=\sum_{\nu=0}^{\infty} \sum_{m\in\zn}{\!}^{\nu, \Omega}\, \lambda_{\nu,m}\,a_{\nu,m}, \quad\mbox{convergence being in } {\mathcal D}'(\Omega) ,
      \eeq
     in the sense of Definition \ref{atomic_domains}, where $(a_{\nu,m})_{\nu \in\nat_0,m\in\zn}$ are non-smooth interior $[K,L]$-atoms, or non-smooth boundary $[K,0]$-atoms according to Definition \ref{def_non-smooth_atoms} and $\lambda\in b_{\p, \q}^{\bm{w}}(\Omega)$. Moreover,
     \beq
	\|f\mid B_{\p, \q}^{\bm{w}}(\Omega)\| \sim \inf \|\lambda \mid b_{\p, \q}^{\bm{w}}(\Omega)\|, \nonumber
     \eeq
     where the infimum is taken over all possible representations of $f$.
     \item[{\upshape (ii)\hfill}] Let $p^+,q^+<\infty$. Then
        \beq Atom(F_{\p, \q}^{\bm{w}})^n \supset R(n).
        \eeq

      Let $K, L\geq 0$ with $K> \alpha_2$ and $L> \sigma_{p,q} - \alpha_1$. Then $f \in \mathcal{D}'(\rn)$ belongs to $F_{\p, \q}^{\bm{w}}(\Omega)$ if and only if it can be represented as
      \beq
	f=\sum_{\nu=0}^{\infty} \sum_{m\in\zn}{\!}^{\nu, \Omega}\, \lambda_{\nu,m}\,a_{\nu,m}, \quad\mbox{convergence being in } {\mathcal D}'(\Omega) ,
      \eeq
     in the sense of Definition \ref{atomic_domains}, where $(a_{\nu,m})_{\nu \in\nat_0,m\in\zn}$ are non-smooth interior $[K,L]$-atoms, or non-smooth boundary $[K,0]$-atoms according to Definition \ref{def_non-smooth_atoms} and $\lambda\in f_{\p, \q}^{\bm{w}}(\Omega)$. Moreover,
     \beq
	\|f\mid F_{\p, \q}^{\bm{w}}(\Omega)\| \sim \inf \|\lambda \mid f_{\p, \q}^{\bm{w}}(\Omega)\|,\nonumber
     \eeq
     where the infimum is taken over all possible representations of $f$.
  \end{list}
\end{theorem}

\begin{remark}
  \begin{list}{}{\labelwidth1.3em\leftmargin2em}
    \item[{\upshape (a)\hfill}] The proof of this result follows in the same way as the proof of \cite[Theorem 3.6]{TW96}, and for that reason we do not present it here. The most important step is the construction of extended non-smooth atoms on $\rn$ according to Definition \ref{defatoms}. The remaining part of the proof is based on the same techniques as used in the proof of the previous result.
    \item[{\upshape (b)\hfill}] Note that the construction in \cite{TW96} ensures moment conditions of the type
    \beq
       \int x^{\beta} a(x) \, dx = 0, \quad \mbox{for } |\beta|<L \mbox{ and } \nu \geq 1. \nonumber
    \eeq
    Since these conditions are more restrictive than the conditions we use here (as stated in Remark \ref{remark_nonsmoothatoms}), the new atoms on $\rn$ are also non-smooth atoms according to Definition \ref{defatoms}.
  \end{list}
\end{remark}
%

\end{document}